\documentclass[12pt,a4paper,reqno]{amsart}
\usepackage{amsmath,latexsym,epsfig,graphics,amsfonts,amssymb}

\newtheorem{theorem}{Theorem}[section]
\newtheorem{lemma}[theorem]{Lemma}

\theoremstyle{definition}

\theoremstyle{remark}

\theoremstyle{plain}

\numberwithin{equation}{section}

\newcommand{\Ho}{\ensuremath{\textup{H}(\mathbb{S})} } 
\newcommand{\Hop}{\ensuremath{\textup{H}^+(\mathbb{S})} }
\newcommand{\Hon}{\ensuremath{\textup{H}^-(\mathbb{S})} }
\newcommand{\So}{\ensuremath{\mathbb{S}} }

\title[Reversibility of circle homeomorphisms]{Reversibility in the group of homeomorphisms of the circle}

\author{Nick Gill, Anthony G. O'Farrell \and Ian Short}

\thanks{The second and third authors were supported by Science Foundation Ireland, grant 05/RFP/MAT0003.}

\begin{document}

\begin{abstract}
We identify those elements of the homeomorphism group of the circle that can be expressed as a composite of two involutions.
\end{abstract}
\maketitle
\section{Introduction}\label{S: intro}

We describe an element $g$ of a group $G$ as \emph{reversible} in $G$ if it is conjugate in $G$ to its own inverse. We say that $g$ is \emph{strongly reversible} in $G$ if there exists an involution $\tau$ in $G$ such that $\tau g \tau = g^{-1}$. This language has developed from the theory of finite groups, where the terms \emph{real} and \emph{strongly real} replace \emph{reversible} and \emph{strongly reversible}. (The word \emph{real} is used because an element $g$ of a finite group is reversible if and only if each irreducible character of $G$ takes a real value when applied to $g$.) Notice that $g$ is strongly reversible if and only if it can be expressed as a composite of two involutions. The strongly reversible elements of the homeomorphism group of the real line were determined by Jarczyk and Young; see \cite{Ja02a,OF04,Yo94}. The purpose of this paper is to determine the strongly reversible maps in the group of homeomorphisms of the circle. 

Let \So denote the unit circle in $\mathbb{R}^2$ centred on the origin. Denote by \Ho the group of homeomorphisms of $\So$. There is a subgroup \Hop of \Ho consisting of orientation preserving homeomorphisms. The subgroup \Hop has a single distinct coset \Hon in \Ho which consists of orientation reversing homeomorphisms. We classify the strongly reversible maps in both groups \Hop and $\Ho$. A classification of the reversible elements of \Hop and \Ho can be extracted from a conjugacy classification in these two groups. We describe the conjugacy classes of \Hop and \Ho in \S\ref{S: conjugacy}, and comment briefly on reversibility.

For points $a$ and $b$ in $\So$, we write $(a,b)$ to indicate the open anticlockwise interval from $a$ to $b$ in $\mathbb{S}$. Let $[a,b]$ denote the closure of $(a,b)$. For a proper open interval $I$ in $\mathbb{S}$, we say that \emph{$u<v$ in $I$} if $(u,v)\subset I$. To classify the strongly reversible maps in \Hop we need the notion of the \emph{signature} of an orientation preserving homeomorphism which has a fixed point. If $f$ is such a homeomorphism, then each point $x$ in $\mathbb{S}$ is either a fixed point of $f$ or else it lies in an open interval component $I$ in the complement of the fixed point set of $f$. The signature $\Delta_f$ of $f$ is the function from $\So$ to $\{-1,0,1\}$ given by the equation
\[
\Delta_f(x) = 
\begin{cases}
1  & \text{if $x< f(x)$ in $I$},\\
0  & \text{if $f(x)=x$},\\
-1 & \text{if $f(x)< x$ in $I$}.
\end{cases}
\]
The strongly reversible maps in \Hop are classified according to the next theorem, proven in \S\ref{S: two+}.

\begin{theorem}\label{T: two+}
An element $f$ of \Hop is strongly reversible if and only if either it is an involution or else it has a fixed point and there is an orientation preserving homeomorphism $h$ of rotation number $\frac12$   such that $\Delta_f = -\Delta_f\circ h$.
\end{theorem}

Elements in \Hop that cannot be expressed as a composite of two involutions (elements that are not strongly reversible) can nevertheless be expressed as a composite of three involutions. The next theorem is proven in \S\ref{S: three+}.

\begin{theorem}\label{T: three+}
Each member of \Hop can be expressed as a composite of three orientation preserving involutions.
\end{theorem}

We move on to describe the strongly reversible elements in the larger group $\Ho$. There are orientation preserving homeomorphisms of \So which are strongly reversible in $\Ho$, but not strongly reversible in $\Hop$. Before we state our theorem on strong reversibility in $\Ho$, we introduce some notation which is explained in more detail in \S\ref{S: conjugacy}. For a  homeomorphism $f$, the \emph{degree} of $f$, denoted $\text{deg}(f)$, is equal to $1$ if $f$ preserves orientation, and $-1$ if $f$ reverses orientation.  Let $\rho(f)$ denote the rotation number of an orientation preserving homeomorphism $f$. The rotation number is an element of $[0,1)$. If $\rho(f)=0$ then $f$ has a fixed point. If $\rho(f)$ is rational then $f$ has a periodic point, in which case the minimal period of $f$ (the smallest positive integer $n$ such that $f^{n}$ has fixed points) is denoted by $n_f$. If $\rho(f)$ is irrational then we denote the minimal set of $f$, that is, the smallest non-trivial $f$ invariant compact subset of $\mathbb{S}$, by $K_f$. This set is either a perfect and nowhere dense subset of $\mathbb{S}$ (a Cantor set) or else equal to $\mathbb{S}$. In the former case we define $I_f$ to be the set of inaccessible points of $K_f$, and in the latter case we define $I_f$ to be $\mathbb{S}$. There is a continuous surjective map $w_f:\So\rightarrow\So$ of degree $1$ with the properties: (i) $w_f$ maps $I_f$ homeomorphically onto $w_f(I_f)$; (ii) $w_f$ maps each closed interval component in the complement of $I_f$ to a point; (iii)  $w_ff=R_\theta w_f$, where $R_\theta$ is the anticlockwise rotation by $\theta=2\pi\rho(f)$.  The next theorem is proven in \S\ref{S: two}.

\begin{theorem}\label{T: two}
Let $f$ be an orientation preserving member of $\Ho$. Either
\begin{enumerate}
\item[(i)] $\rho(f)=0$, in which case $f$ is strongly reversible if and only if there is a homeomorphism $h$ such that $\Delta_f = -\textup{deg}(h)\cdot \Delta_f\circ h$, and either $h$ preserves orientation and has rotation number $\frac12$ or else reverses orientation;
\item[(ii)] $\rho(f)$ is non-zero and rational, in which case $f$ is strongly reversible if and only if $f^{n_f}$ is strongly reversible by an orientation reversing involution;
\item[(iii)] $\rho(f)$ is irrational, in which case $f$ is strongly reversible if and only if $w_f(I_f)$ has a reflectional symmetry.
\end{enumerate}
\end{theorem}

Using Theorem~\ref{T: two} we also show that an orientation preserving circle homeomorphism is strongly reversible by an orientation \emph{reversing} involution if and only if it is reversible by an orientation \emph{reversing} homeomorphism. 

It remains to state a result on strong reversibility of orientation reversing homeomorphisms. Each orientation reversing homeomorphism has exactly two fixed points. The next theorem is proven in \S\ref{S: two-}.

\begin{theorem}\label{T: two-}
An orientation reversing homeomorphism $f$ is strongly reversible if and only if  there is an orientation reversing homeomorphism $h$ that interchanges the pair of fixed points of $f$ and satisfies $h f^2h^{-1} = f^{-2}$.
\end{theorem}

Fine and Schweigert proved in \cite[Theorem 25]{FiSc55} that each member of \Ho can be expressed as a composite of three involutions  (this result follows quickly from Theorems \ref{T: three+} and \ref{T: two}).

We now describe the structure of this paper. Section \ref{S: conjugacy}  contains unoriginal material; it consists of a brief review of a conjugacy classification in \Hop and $\Ho$. All subsequent sections contain new results. Sections \ref{S: two+} and \ref{S: three+} are about \Hop and sections \ref{S: two} and \ref{S: two-} are about $\Ho$.

\section{Conjugacy classification}\label{S: conjugacy}

Two conjugacy invariants which can be used to determine the conjugacy classes in \Ho are the rotation number and signature, introduced in \S\ref{S: intro}. We describe here only those properties of these two quantities that we will use. For more information on rotation numbers, see \cite{Gh01}; for more information on signatures, see \cite{FiSc55}.  

For an orientation preserving homeomorphism $f$, choose any point $x$ in $\mathbb{S}$, and let $\theta_n$ be the angle in $[0,2\pi)$ measured anticlockwise between $f^{n-1}(x)$ and $f^n(x)$. The \emph{rotation number} of $f$, denoted $\rho(f)$, is the unique number in $[0,1)$ such that the expression
\[
(\theta_1+\dots+\theta_n) - 2\pi n\rho(f)
\]
is bounded for all $n$. The quantity $\rho(f)$ is independent of $x$. The rotation number is invariant under conjugation in $\Hop$, and  $\rho(f^n)=n\rho(f) \pmod {1}$, for each integer $n$. 

A straightforward consequence of the definition of $\rho(f)$ is that $\rho(f)=0$ if and only if $f$ has a fixed point. In this case, $\mathbb{S}$ can be partitioned into a closed set $\text{fix}(f)$, consisting of fixed points of $f$, and a countable collection of open intervals on each of which $f$ is free of fixed points. Now suppose that $\rho(f)=p/q$, where $p$ and $q$ are coprime positive integers. Then $f$ has periodic points, that is, there is a positive integer $n$ for which $f^n$ has fixed points. The smallest such $n$, denoted $n_f$, is the \emph{minimal period} of $f$ and is equal to $q$. 

The remaining possibility is that $\rho(f)$ is irrational. In this case we define $K_f$ to be the unique minimal set in the poset consisting of $f$ invariant compact subsets of $\mathbb{S}$ ordered by inclusion. We describe $K_f$ as the \emph{minimal set} of $f$. Either $K_f=\mathbb{S}$ or else $K_f$ is a perfect subset of $\mathbb{S}$ with empty interior---a Cantor set. In the latter case
there is a sequence of open intervals $(a_i,b_i)$, for $i=1,2,\dotsc$, such that $[a_i,b_i]\cap [a_j,b_j]=\emptyset$ when $i\neq j$,  and $K_f$ is the complement of $\bigcup_{i=1}^\infty (a_i,b_i)$. The set $I_f$ of \emph{inaccessible points} of $K_f$ is the complement of 
$\bigcup_{i=1}^\infty [a_i,b_i]$. If $K_f=\mathbb{S}$ then we define $I_f=\mathbb{S}$. There is a continuous surjective map $w_f$ of $\mathbb{S}$ of degree $1$ such that $w_ff=R_\theta w_f$, where $\theta =2\pi\rho(f)$. The map $w_f$ is a homeomorphism when restricted to $I_f$, and it maps each interval $[a_i,b_i]$ to a single point. The map $w_f$
is unique up to post composition by rotations. 

Now suppose that $f$ is an orientation preserving homeomorphism which has a fixed point. The signature of $f$ was defined in \S\ref{S: intro}. The signature $\Delta_f$ takes the value $0$ on $\text{fix}(f)$, and elsewhere it takes either the value $-1$ or the value $1$. Useful properties of the signature are encapsulated in the next elementary lemma.

\begin{lemma}\label{L: handy}
If $f$ is a member of $\Hop$ with a fixed point, and $h$ is a member of $\Ho$, then 
\begin{enumerate}
\item[(i)] $\Delta_{hfh^{-1}}=\textup{deg}(h)\cdot\Delta_f\circ h^{-1}$,
\item[(ii)] $\Delta_{f^{-1}}=-\Delta_{f}$.
\end{enumerate}
\end{lemma}

We are now in a position to state criteria which determine whether two circle homeomorphisms are conjugate. The results are stated in such a way that one can deduce from them when two orientation preserving circle homeomorphisms are conjugate in each of the groups \Hop and $\Ho$. The result on irrational rotation numbers follows from \cite[Theorem 2.3]{Ma70}. The result on orientation preserving homeomorphisms with fixed points is similar to \cite[Theorem 10]{FiSc55}. The other two theorems are well-known; they can both be proven directly. Recall that the degree of a circle homeomorphism is $1$ if the map preserves orientation, and $-1$ if it reverses orientation.

\begin{theorem}\label{T: conjFix}
Two orientation preserving circle homeomorphisms $f$ and $g$, each of which has a fixed point, are conjugate by a homeomorphism of degree $\epsilon$ if and only if there is a homeomorphism $h$ of degree $\epsilon$ such that $\Delta_g = \epsilon\Delta_f\circ h$.
\end{theorem}

Note in particular that a map $g$ in $\Hop$ with fixed points is conjugate in $\Hop$ to all of its powers.

\begin{theorem}\label{T: conjPer}
Two orientation preserving circle homeomorphisms $f$ and $g$, both of which have the same non-zero rational rotation number, are conjugate by a homeomorphism of degree $\epsilon$ if and only if $f^{n_f}$ is conjugate to $g^{n_g}$ by a homeomorphism of degree $\epsilon$.
\end{theorem}

Since $f$ and $g$ have the same non-zero rational rotation number, the integers $n_f$ and $n_g$ in Theorem~\ref{T: conjPer} are equal. 

Recall that an orthogonal map of the circle of degree $1$ is a rotation, and an orthogonal map of the circle of degree $-1$ is a reflection in a line through the origin.

\begin{theorem}\label{T: conjIrrat}
Two orientation preserving circle homeomorphisms $f$ and $g$, both of which have the same irrational rotation number, are conjugate by a homeomorphism of degree $\epsilon$ if and only if there is an orthogonal map of degree $\epsilon$ that maps $w_f(I_f)$ to $w_g(I_g)$.
\end{theorem}
 
It remains to consider conjugacy between orientation reversing maps.

\begin{theorem}\label{T: conj-}
Two orientation reversing circle homeomorphisms $f$ and $g$ are conjugate in $\Ho$ if and only if $f^2$ and $g^2$ are conjugate in $\Ho$ by a homeomorphism that maps the pair of fixed points of $f$ to the pair of fixed points of $g$.
\end{theorem}

It follows from Theorems~\ref{T: conjPer} and \ref{T: conj-} that all non-trivial involutions in $\mathcal{H}^+$ are conjugate, and all orientation reversing involutions in $\mathcal{H}$ are conjugate. (These statements can easily be seen directly.)

We briefly remark on the reversible elements in \Hop and $\Ho$. 
Suppose that $f$ and $h$ are members of \Hop such that  $hfh^{-1}=f^{-1}$. Then
\[
\rho(f^{-1})=\rho(hfh^{-1}) = \rho(f)\pmod {1}. 
\]
But $\rho(f^{-1})=-\rho(f)\pmod {1}$, hence $\rho(f)$ is equal to either $0$ or $\frac{1}{2}$. One can construct examples of reversible elements in \Hop with either of these two rotation numbers (there are examples at the end of \S\ref{S: two}). On the other hand, if $h$ reverses orientation and still $hfh^{-1}=f^{-1}$, then
\[
\rho(hfh^{-1}) = -\rho(f)\pmod {1}, 
\]
so, in this case, the rotation number tells us nothing about reversibility. Notice that if we compare Theorems~\ref{T: conjFix}, \ref{T: conjPer}, and \ref{T: conjIrrat} with Theorem~\ref{T: two}, we see that that an orientation preserving map $f$ is reversible by an orientation \emph{reversing} map if and only if $f$ is strongly reversible by an orientation \emph{reversing} involution. There are, however, orientation preserving homeomorphisms that are not strongly reversible in $\mathcal{H}$, but are nevertheless reversible by orientation \emph{preserving} maps; one example is given at the end of \S\ref{S: two}.

\section{Proof of Theorem~\ref{T: two+}}\label{S: two+}

 The following two theorems deal with strong reversibility in \Hop for the two rotation numbers $0$ and $\frac{1}{2}$ separately.

A result similar to Theorem~\ref{T: rot0}, below, has been proven by Jarzcyk \cite{Ja02a} and Young \cite{Yo94} for homeomorphisms of the real line. We use an elementary lemma in the proof of Theorem~\ref{T: rot0}.

\begin{lemma}\label{L: elementary}
If $f$ and $g$ are orientation preserving homeomorphisms such that $\Delta_f=\Delta_g$ then there is an orientation preserving homeomorphism $k$ that fixes each of the fixed points of $f$ and $g$ such that $kfk^{-1}=g$.
\end{lemma}
\begin{proof}
We define a homeomorphism $k$ as follows. On each fixed point $x$ of $f$ and $g$, define $k(x)=x$. On each open interval component $(a,b)$ of $\mathbb{S}\setminus\text{fix}(f)$, the signature function takes either the value $1$ for both functions $f$ and $g$, or else it takes the value $-1$ for both functions. In either case we can choose an orientation preserving homeomorphism $k_0$ of $(a,b)$ such that $k_0 fk_0^{-1}=g$ for $x\in(a,b)$. We then define $k(x)=k_0(x)$ for $x\in(a,b)$. We have constructed the required function $k$. Of course, the existence of a conjugation between $f$ and $g$ follows from Theorem~\ref{T: conjFix}, but we also needed the property of $k$ that it fixes each element of $\text{fix}(f)$.  
\end{proof}

\begin{theorem}\label{T: rot0}
An element  $f$ of $\Hop$ with a fixed point is strongly reversible in \Hop if and only if there is a homeomorphism $h$ in \Hop with rotation number $\frac12$  such that $\Delta_f=-\Delta_f\circ h$.
\end{theorem}
\begin{proof}
If $\sigma f\sigma =f^{-1}$ for an orientation preserving involution $\sigma$, then $\Delta_f = -\Delta_f\circ\sigma$, by Lemma~\ref{L: handy} (i). Conversely, suppose that there is a homeomorphism $h$ in \Hop with rotation number $\frac12$ such that $\Delta_f=-\Delta_f\circ h$.
By Lemma~\ref{L: handy}, $\Delta_{h^{-1} fh}=\Delta_{f^{-1}}$. Using Lemma~\ref{L: elementary} we can construct a map $k$ in $\Hop$ that fixes each fixed point of $f$, and satisfies $k^{-1}h^{-1} fh k = f^{-1}$. 

Now choose a fixed point $p$ of $f$. Then $\Delta_{f}(h(p))=-\Delta_{f}(p)=0$, so $h(p)$ is a fixed point of $f$.  The points $p$ and $h(p)$ are distinct because $h$ has no fixed points. Define
\[
\mu(x)=
\begin{cases}
hk(x) & \text{if $x\in[p,h(p)]$},\\
k^{-1} h^{-1}(x) & \text{if $x\in[h(p),p]$}.
\end{cases}
\] 
One can check that $\mu$ is an involution in \Hop and $\mu f\mu=f^{-1}$.
\end{proof}

We move on to orientation preserving homeomorphisms with rotation number $\frac{1}{2}$.

\begin{theorem}\label{T: rot1/2}
An element of $\Hop$ with rotation number $\frac{1}{2}$ is strongly reversible in $\Hop$ if and only if it is an involution.
\end{theorem}
\begin{proof}
All involutions are strongly reversible by the identity map.  Conversely, let $f$ be a homeomorphism with rotation number $\frac12$, and let $\sigma$ be an involution in  $\Hop$ such that $\sigma f\sigma =f^{-1}$. Choose an element $x$ of $\text{fix}(f^2)$. Then $f(x)$ is also an element of $\text{fix}(f^2)$, and by interchanging $x$ and $f(x)$ if necessary, we may assume that $\sigma(x) \in (x,f(x)]$. Suppose that $\sigma(x)\neq f(x)$. Since $\sigma$ maps $(\sigma(x),x)$ onto $(x,\sigma(x))$, we have that $\sigma f(x)\in (x,\sigma(x))$. Likewise, $f$ maps $(x,f(x))$ onto $(f(x),x)$, therefore $f\sigma(x) \in (f(x),x)$. However, $f\sigma(x)= \sigma f ^{-1}(x)=\sigma f(x)$, and yet $(x,\sigma(x)) \cap (f(x),x)=\emptyset$. This is a contradiction, therefore $\sigma(x)=f(x)$. This means that $\sigma f$ is an orientation preserving involution which fixes $x$; hence it is the identity map. Therefore $f=\sigma$, as required. 
\end{proof}

In contrast to Theorem~\ref{T: rot1/2} there are reversible homeomorphisms with rotation number $\frac{1}{2}$ that are not strongly reversible. An example is given at the end of \S\ref{S: two}.

We have all the ingredients for a proof of Theorem~\ref{T: two+}.

\begin{proof}[Proof of Theorem~\ref{T: two+}]
If $f$ is a strongly reversible member of \Hop then it is reversible, so it must have rotation number equal to either $0$ or $\frac{1}{2}$. If it has rotation number $0$ then there is an orientation preserving homeomorphism $h$ with rotation number $\frac12$ $\sigma$  such that $\Delta_f=-\Delta_f\circ h$, by Theorem~\ref{T: rot0}. If $f$ has rotation number $\frac{1}{2}$ then it is an involution, by Theorem~\ref{T: rot1/2}. The converse implication follows immediately from Theorem~\ref{T: rot0}.
\end{proof}

\section{Proof of Theorem~\ref{T: three+}}\label{S: three+}

\begin{proof}[Proof of Theorem~\ref{T: three+}]
Choose an element $f$ of \Hop that is not an involution. There exists a point $x$ in $\mathbb{S}$ such that $x$, $f(x)$, and $f^2(x)$ are three distinct points. By replacing $f$ with $f^{-1}$ if necessary we can assume that $x$, $f(x)$, and $f^2(x)$ occur in that order anticlockwise around $\mathbb{S}$. Notice that $f^{-1}(x)$ lies in $(f(x),x)$. Choose a point $y$ in $(x,f(x))$ that is sufficiently close to $f(x)$ that $f^{-1}(y)>f^2(x)$ in $(f(x),x)$. We construct an orientation preserving homeomorphism $g$ from $[x,f(x)]$ to $[f(x),x]$ such that $g(y)=f(y)$, $g(t)<\text{min}(f(t),f^{-1}(t))$ in $(x,y)$,	 and $f(t)<g(t)<f^{-1}(t)$ in $(y,f(x))$. A graph of such a function is shown in Figure~\ref{F: g}.
\begin{figure}[ht]
\centering
\includegraphics[scale=0.7]{g.eps}
\caption{}
\label{F: g}
\end{figure}
Define an involution $\sigma$ in $\Hop$ by the equation
\[
\sigma(t)=
\begin{cases}
g(t) & \text{if $t\in[x,f(x)]$},\\
g^{-1}(t) & \text{if $t\in[f(x),x]$}.
\end{cases}
\]
Let us determine the fixed points of $\sigma f$.  For a point $t$ in $[x,f(x)]$ we have that $\sigma f(t)=t$ if and only if $f(t)=g(t)$. This means that either $t=x$ or $t=y$. For a point $t$ in $(f(x),x)$ we have that $\sigma f(t)=t$ if and only if $\sigma f g(w)= g(w)$, where $w=g^{-1}(t)$ is a point in $(x,f(x))$. Therefore $g(w)=f^{-1}(w)$. This equation has no solutions in $(x,f(x))$, hence $\sigma f$ has no fixed points in $(f(x),x)$. It is straightforward to obtain the direction of flow on the complement of $\{x,y\}$ and we find that
\[
\Delta_{\sigma f}(t)=
\begin{cases}
0  & \text{if $t=x,y$},\\
1  & \text{if $t\in(x,y)$},\\
-1 & \text{if $t\in(y,x)$}.
\end{cases}
\]
By Theorem~\ref{T: rot0}, $\sigma f$ is expressible as a composite of two involutions in $\Hop$. Therefore $f$ is expressible as a composite of three involutions in $\Hop$. 
\end{proof}

A simple corollary of Theorem~\ref{T: three+} is that \Hop is uniformly perfect, meaning that there is a positive integer $N$ such that each element of \Hop can be expressed as a composite of $N$ or fewer commutators. Since it is easy to express the rotation by $\pi$ as a commutator, and each involution in \Hop is conjugate to the rotation by $\pi$, it follows from Theorem~\ref{T: three+} that \Hop is uniformly perfect with $N=3$. In fact, Eisenbud, Hirsch, and Neumann \cite{EiHiNe81} proved that \Hop is uniformly perfect with $N=1$.

\section{Proof of Theorem~\ref{T: two}}\label{S: two}

For the remainder of this document we work in the full group of homeomorphisms of the circle. In \S\ref{S: conjugacy} we showed that all reversible maps in \Hop have rotation number either $0$ or $\frac{1}{2}$. This is not the case in \Ho because if $h$ is an orientation reversing map, and $f$  an orientation preserving map, then $\rho(hfh^{-1})=-\rho(f)$. Since also $\rho(f^{-1})=-\rho(f)$, the rotation number tells us nothing about reversibility by orientation reversing homeomorphisms.  In fact, since all rotations are strongly reversible in $\Ho$ by reflections, there are strongly reversible maps in $\Ho$ with any given rotation number. 

We divide our analysis of strongly reversible maps in \Ho between three cases corresponding to when the rotation number is $0$, rational, or irrational. The first case is Theorem~\ref{T: two} (i). We need a preliminary lemma.

\begin{lemma}\label{L: preliminary}
If $f$ is a fixed point free homeomorphism of an open proper arc $A$ in $\mathbb{S}$, then $f$ is conjugate to $f^{-1}$ on $A$ by an orientation reversing involution of $A$.
\end{lemma}
\begin{proof}
The situation is topologically equivalent to the situation when $f$ is a fixed point free homeomorphism of the real line. Such maps $f$ are conjugate to non-trivial translations, and translations are reversible by the orientation reversing involution $x\mapsto -x$.
\end{proof}

\begin{proof}[Proof of Theorem~\ref{T: two} (i)]
If $\tau f\tau =f^{-1}$ for an involution $\tau$, then 
 $\Delta_f = -\text{deg}(\tau)\cdot \Delta_f\circ\tau$, by Lemma~\ref{L: handy}. For the converse, we are given a homeomorphism $h$ that satisfies $\Delta_f=-\text{deg}(h)\cdot\Delta_f\circ h$. Either $h$ preserves orientation and satisfies $\rho(h)=\frac12$, in which case the result follows from Theorem~\ref{T: rot0}, or else $h$ reverses orientation. 

In the latter case, by Lemmas~\ref{L: handy} and \ref{L: elementary} there is an orientation preserving homeomorphism $k$ that fixes the fixed points of $f$, and satisfies $k^{-1}h^{-1} fh k = f^{-1}$. Define $s=k^{-1}h^{-1}$. Let $s$ have fixed points $p$ and $q$. Let $a$ denote the point in $\text{fix}(f)$ that is clockwise from $p$, and closest to $p$. Possibly $a=p$. Define $b$ to be the point in $\text{fix}(f)$ that is anticlockwise from $p$ and closest to $p$. Let $I=(a,b)$. Similarly we define an interval $J$ about $q$. If we ignore the trivial case in which $\text{fix}(f)$ has only one component then $I$ and $J$ only intersect, if at all, in their end-points. Now, $f$ fixes $I$ so we can, by Lemma~\ref{L: preliminary}, choose an orientation reversing involution $\tau_I$ of $I$ such that $\tau_I f\tau_I(x)=f^{-1}(x)$ for $x\in I$. Similarly we define $\tau_J$. From the equation $sfs^{-1}=f^{-1}$ we deduce that $s$ fixes $I$ and $J$. Hence we can define
\begin{equation}\label{E: invert}
\mu(x) = 
\begin{cases}
\tau_I(x) & \text{if $x\in I$},\\ 
\tau_J(x) & \text{if $x\in J$},\\
s(x) & \text{if $x\in [p,q]\setminus (I\cup J)$},\\
s^{-1}(x) & \text{if $x\in [q,p]\setminus (I\cup J)$}. 
\end{cases}
\end{equation}
One can check that $s$ is an orientation reversing involution that satisfies $sfs=f^{-1}$.
\end{proof}

To prove Theorem~\ref{T: two} (ii) we use a lemma that enables us to deal with rational rotation numbers of the form $1/n$, rather than $m/n$. Recall that $n_f$ denotes the minimal period of $f$.

\begin{lemma}\label{L: simplify}
Let $f$ be an element of $\Hop$ with a periodic point. If $f^d$ is strongly reversible for an integer $d$ in $\{1,2,\dots,n_f-1\}$ that is coprime to $n_f$, then $f$ is strongly reversible. 
\end{lemma}
\begin{proof}
There exist integers $u$ and $t$ such that $dt= 1+un_f$. Let $q$ be the positive integer between $0$ and $n_f$, and coprime to $n_f$, such that $\rho(f)=q/n_f$. Observe that 
\[
\rho(f^{dt})=(1+un_f)\rho(f)=\rho(f)+uq = \rho(f) \pmod {1}.
\]
Recall that a map $g$ in $\Hop$ with fixed points is conjugate in \Hop to all its powers. Let $g=f^{n_f}$. Then $g$ is conjugate to $g^{dt}$. In other words, $f^{n_f}$ is conjugate to $(f^{dt})^{n_f}$. Apply Theorem~\ref{T: conjPer} to the maps $f$ and $f^{dt}$ to see that these two maps are conjugate. The second map $f^{dt}$ is strongly reversible, because $f^d$ is strongly reversible. Conjugacy preserves strong reversibility, therefore $f$ is also strongly reversible.
\end{proof}

We first prove a special case of Theorem~\ref{T: two} (ii).

\begin{lemma}\label{L: reflect}
Let $f$ be an orientation preserving homeomorphism of $\mathbb{S}$ with rotation number $1/n$, for a positive integer $n$. Suppose that there is an open interval $J$ such that the intervals $J, f(J),\dots,f^{n-1}(J)$ are pairwise disjoint, and such that $\textup{fix}(f^n)$ is the complement of $J\cup f(J)\cup\dots\cup f^{n-1}(J)$. Then $f$ is strongly reversible by an orientation reversing involution that maps $f^k(J)$ to $f^{-k+1}(J)$ for each integer $k$.
\end{lemma}
\begin{proof}
Let $R$ denote an anticlockwise rotation by $2\pi/n$. After conjugating $f$ suitably, the function $f$ and interval $J$ may be adjusted so that $f^k(J)=R^k(J)$ for all integers $k$.  Let $\tau$ denote reflection in a line $\ell$ through the origin that bisects $J$. Thus $\tau$ fixes $J$. Let $\sigma$ denote reflection in a line through the origin that is $\pi/n$ anticlockwise from $\ell$. Then $R=\sigma\tau$. Orient $J$ in an anticlockwise sense, and choose an increasing homeomorphism $\phi$ of $J$, without fixed points, such that $\tau\phi\tau=\phi^{-1}$. This is possible because we can, by conjugation, consider $J$ to be the real line and consider $\tau$ still to be a reflection, in which case $\phi$ can be chosen to be a translation. Now define an orientation preserving circle homeomorphism $g$ to satisfy $g(x)=f(x)$, for $x\in\text{fix}(f^n)$, and $g(x)=R^{k+1}\phi R^{-k}(x)$, for $x\in R^k(J)$. This means that $g^n$ has the same signature on all of the intervals $R^k(J)$ (either all $-1$ or all $1$). By Theorem~\ref{T: conjPer}, $g$ is conjugate to $f$. Also, one can check that $\sigma g\sigma =g^{-1}$ and $\sigma (g^k(J))=g^{-k+1}(J)$ for each integer $k$. Since strong reversibility is preserved under conjugation, the result follows.  
\end{proof}

\begin{proof}[Proof of Theorem~\ref{T: two} (ii)]
 If $f$ is strongly reversible by an involution $\tau$ then $f^{n_f}$ is also strongly reversible by $\tau$. If $\tau$ is orientation preserving then, by Theorem \ref{T: two+}, $f$ is an involution. Therefore $f^{n_f}$ is the identity, and as such it is reversible by any orientation reversing involution.

Conversely, suppose  that there is an orientation reversing involution $\tau$ such that $\tau f^{n_f}\tau = f^{-n_f}$. Let $p$ be a fixed point of $f^{n_f}$. There is an integer $d$ that is coprime to $n_f$ such that the distinct points $p, f^d(p),\dots,f^{(n_f-1)d}(p)$ occur in that order  anticlockwise around $\mathbb{S}$. The function $(f^d)^{n_f}$ is strongly reversible because it equals $(f^{n_f})^d$. If we can deduce that $f^d$ is strongly reversible then it follows from Lemma~\ref{L: simplify} that $f$ is strongly reversible. In other words, it is sufficient to prove the theorem when the  points $p,f(p),\dots,f^{{n_f}-1}(p)$ occur in that order around $\mathbb{S}$. 

The map $\tau f$ has a fixed point $q$ which, by replacing $\tau$ with $f^k\tau f^{-k}$ for an appropriate integer $k$, we can assume that it lies in the interval $(f^{-1}(p),p]$. Either $q$ is a fixed point of $f^{n_f}$ (that is, a periodic point of $f$) or it is not. In the former case let $I=[q,f(q)]$ and define, for each integer $k$,
\begin{equation}\label{E: mu}
\mu(x)= f^k\tau f^k(x),\quad x\in f^{-k}(I).
\end{equation}
This is a well defined homeomophism because $f^{n_f}\tau f^{n_f}=f^{-n_f}$. One can check that $\mu$ is an involution and satisfies $\mu f\mu =f^{-1}$. In Figure~\ref{F: qIterates} the action of $\mu$ on certain $f$ iterates of $q$ is shown in the case $n_f=2m$.
\begin{figure}[ht]
\centering
\includegraphics[scale=0.7]{qIterates.eps}
\caption{}
\label{F: qIterates}
\end{figure}
If $q$ is not a fixed point of $f^{n_f}$ then it lies in a unique component $J=(s,t)$ in the complement of $\text{fix}(f^{n_f})$. The interval $J$ is contained in $(f^{-1}(p),p)$, and both $s$ and $t$ lie in $\text{fix}(f^{n_f})$. From the equation $\tau f^{n_f}\tau = f^{-n_f}$, we can deduce that $\tau$ maps $J$ to another open interval component in $\mathbb{S}\setminus\text{fix}(f^{n_f})$. Also, $f(J)$ is a component of $\mathbb{S}\setminus\text{fix}(f^{n_f})$. But $\tau(J)$ and $f(J)$ both contain the point $\tau(q)=f(q)$; therefore $\tau(J)=f(J)$. This means that $\tau(s)=f(t)$ and $\tau(t)=f(s)$.

We define an involution $\mu$ in a similar fashion to \eqref{E: mu}, but this time we do not, yet, define $\mu$ on the intervals $J,f(J),\dots, f^{n_f-1}(J)$. Specifically,  let $I=[t,f(s)]$ and define, for each integer $k$, 
\begin{equation*}\label{E: mu1}
\mu(x)= f^k\tau f^k(x),\quad x\in f^{-k}(I).
\end{equation*}
We can extend $\mu$ to $J,f(J),\dots, f^{n_f-1}(J)$ using Lemma~\ref{L: reflect}. The now fully defined map $\mu $ is an orientation reversing involutive homeomorphism of $\mathbb{S}$ that satisfies $\mu f\mu=f^{-1}$.
\end{proof}

To prove Theorem~\ref{T: two} (iii) we need a lemma. We prove the lemma explicitly, although it can be deduced quickly from Lemma~\ref{L: reflect}.

\begin{lemma}\label{L: involution}
Let $J$ and $J'$ be two disjoint non-trivial closed intervals in $\mathbb{S}$, and let $g$ be an orientation preserving homeomorphism from $J$ to $J'$. Then there exists an orientation reversing homeomorphism $\gamma$ from $J$ to $J'$ such that $\gamma g^{-1} \gamma =g$.
\end{lemma}
\begin{proof}
Let $a$ and $b$ be points such that $J=[a,b]$. Choose a point $q$ in $(a,b)$. Choose an orientation reversing homeomorphism $\alpha$ from $[a,q]$ to $[g(q),g(b)]$. The map $\gamma$ defined by
\begin{equation*}
\gamma(x) = 
\begin{cases}
\alpha(x) & \text{if $x\in [a,q]$},\\
g\alpha^{-1}g(x) & \text{if $x\in [q,b]$},
\end{cases}
\end{equation*}
has the required properties.
\end{proof}

\begin{proof}[Proof of Theorem~\ref{T: two} (iii)]
If $f$ is strongly reversible (by an orientation reversing involution) then, since $I_{f^{-1}}=I_f$, we see from Theorem~\ref{T: conjIrrat} that $w_f(I_f)$ has a reflectional symmetry. 

Conversely, suppose that there is a reflection $\tau$ in a line through the origin of $\mathbb{R}^2$ that fixes $w_f(I_f)$. We define an involution $\mu$ from $I_f$ to $I_f$ by the equation $\mu(x) = w_f^{-1}\tau w_f(x)$. In this equation, $w_f^{-1}$ is the inverse of the function $w_f:I_f\rightarrow w_f(I_f)$. For $x\in I_f$ we have
\[
\mu f \mu(x) = (w_f^{-1}\tau w_f)(w_f^{-1}R_\theta w_f)(w_f^{-1}\tau w_f)(x) = f^{-1}(x). 
\] 

Recall that $K_f$ is the complement in $\mathbb{S}$ of a countable collection of disjoint open intervals $(a_i,b_i)$, and $I_f$ is the complement in $\mathbb{S}$ of the union of the intervals $[a_i,b_i]$. We can extend the definition of $\mu$ to $K_f$ by defining $\mu(a_i)$ to be the limit of $\mu(x_n)$, where $x_n$ is a sequence in $I_f$ that converges to $a_i$. Similarly for $b_i$. The extended map $\mu$ is a homeomorphism from $K_f$ to itself. Notice that $\mu$ has the property that for each integer $i$ there is an integer $j$ such that $\mu$ interchanges $a_i$ and $b_j$, and  also interchanges $a_j$ and $b_i$.

 We can extend the definition of $\mu$ to the whole of $\mathbb{S}$ by introducing, for each $i$, an orientation reversing homeomorphism $\phi_i:[a_i,b_i]\rightarrow [a_j,b_j]$ (where $\mu(a_i)=b_j$ and $\mu(b_i)=a_j$), and defining $\mu(x)=\phi_i(x)$ for $x\in[a_i,b_i]$. The resulting map $\mu$ will be a homeomorphism. It remains only to show how to choose particular maps $\phi_i$ such that $\mu$ is an involution that satisfies $\mu f\mu = f^{-1}$.

Let $J=(a_i,b_i)$ and let $J'=(a_j,b_j)$. We have two collections   
\begin{equation}\label{E: intervals}
\{\dotsc,f^{-1}(J),J,f(J),f^2(J),\dotsc\},\quad \{\dotsc,f^{-1}(J'),J',f(J'),f^2(J'),\dotsc\},
\end{equation}
each consisting of pairwise disjoint intervals. These two collections either share no common members, or else they coincide. In the first case, choose an arbitrary orientation reversing homeomorphism $\gamma$ from $J$ to $J'$. In the second case, there is an integer $m$ such that $f^m(J)=J'$. Apply Lemma~\ref{L: involution} with $g=f^m$ to deduce the existence of an orientation reversing homeomorphism $\gamma$ from $J$ to $J'$ satisfying $\gamma f^{-m} \gamma =f^{m}$. In each case we define, for each $n\in\mathbb{Z}$,
\begin{equation*}\label{E: kappa}
\mu(x) = 
\begin{cases}
f^n\gamma f^n(x) & \text{if $x\in f^{-n}(J)$},\\
f^n\gamma^{-1} f^n(x) & \text{if $x\in f^{-n}(J')$}.
\end{cases}
\end{equation*}
One can check that $\mu$ is well-defined for points $x$ in one of the intervals of \eqref{E: intervals}, and that $\mu^2(x)=x$ and $\mu f\mu(x)=f^{-1}(x)$. In this manner $\mu$ can be defined on each of the intervals $(a_k,b_k)$. The resulting map is a homeomorphism of $\So$ that is an involution and satsifies $\mu f\mu =f^{-1}$.
\end{proof}

It follows from Theorem~\ref{T: two} and the results on conjugacy in \S\ref{S: conjugacy} that an orientation preserving circle homeomorphism is reversible by an orientation \emph{reversing} involution if and only if it is reversible by an orientation \emph{reversing} homeomorphism. We now sketch the details of an example to show that there are orientation preserving circle homeomorphisms that are reversible by orientation preserving homeomorphisms, but are \emph{not} strongly reversible in $\mathcal{H}$. This means that the concepts of reversibility and strong reversibility are not equivalent in either $\mathcal{H}$ or $\mathcal{H}^+$. 

Let $a$ be the point on $\mathbb{S}$ with co-ordinates $(0,1)$ and let $b$ be the point with  co-ordinates $(0,-1)$. Let $\dotsb<a_{-2}<a_{-1}<a_0<a_1<a_2< \dotsb$ be an infinite sequence of points in $(b,a)$ that accumulates only at $a$ and $b$.  Let $g$  be an orientation preserving homeomorphism from $[b,a]$ to $[b,a]$ that fixes  only  the points $a_i$, $a$, and $b$. We construct a doubly infinite sequence $s$ consisting of $1$s and $-1$s as follows. Let $u$ represent the string of six numbers $1,1,1,-1,-1,1$. Let $-u$ represent the string $-1,-1,-1,1,1,-1$. Then $s$ is given by $\dotsc,u,-u,u,-u,\dotsc$. We say that two doubly infinite sequences $(x_n)$ and $(y_n)$ are equal if and only if there is an integer $k$ such that $x_{n+k}=y_n$ for all $n$. Our sequence $s$ has been constructed such that $s=-s$, and the sequence formed by reversing $s$ is distinct from $s$.

Suppose that $g$ is defined in any fashion on the intervals 
\[
\dotsc,(a_{-2},a_{-1}),(a_{-1},a_{0}),(a_{0},a_{1}),(a_{1},a_{2}),\dotsc
\]
such that the signature of $g$ on these intervals
 is determined by $s$. For  $x\in(a,b)$ we define $g(x)=R_{\pi}gR_{\pi}(x)$, where $R_{\pi}$ is the rotation by $\pi$. A diagram of the homeomorphism $g$ is shown in Figure~\ref{F: counterexample}.
\begin{figure}[ht]
\centering
\includegraphics[scale=0.6]{counterexample.eps}
\caption{}
\label{F: counterexample}
\end{figure}
Since we can embed $g$ in a flow, we can certainly choose a square-root $\sqrt{g}$ of $g$ (which shares the same signature as $g$). Define another circle homeomorphism $f$ by the formula
\begin{equation*}
f(x) = 
\begin{cases}
R_{\pi} \sqrt{g}(x) & \text{if $x\in [b,a]$},\\
\sqrt{g}R_{\pi}(x) & \text{if $x\in (a,b)$}.
\end{cases}
\end{equation*}
This map $f$ has rotation number $\frac12$ and it satisfies $f^2=g$. The map $f$ is not reversible by an orientation preserving involution, by Theorem~\ref{T: rot1/2}, because it is not an involution. Nor is $f$ reversible by an orientation reversing involution; for it were then $g$ would also be reversed by the same orientation reversing involution, and from Theorem~\ref{T: two} one can deduce that this would mean that the sequence $s$ coincides with the reversed sequence $s$. Finally, $f$ is reversible since $g$ is reversible, by Theorem~\ref{T: conjPer}; a conjugation from $g$ to $g^{-1}$ can be constructed that maps $a_i$ to $a_{i+6}$ for each $i$.

\section{Proof of Theorem~\ref{T: two-}}\label{S: two-}

\begin{proof}[Proof of Theorem~\ref{T: two-}]
If $f$ is strongly reversible then it is expressible as a composite of two involutions, one of which, $\tau$,  must reverse orientation. From the equation $\tau f\tau = f^{-1}$ we see that $\tau f^2\tau=f^{-2}$, and that $\tau$ preserves the pair of fixed points of $f$ as a set. If $\tau$ fixes each of the fixed points of $f$ then $\tau f$ is an orientation preserving involution with fixed points. Hence it is the identity map. Therefore $f=\tau$. The alternative is that $\tau$ interchanges the fixed points of $f$, which is the condition stated in Theorem~\ref{T: two-}.

For the converse, suppose that $a$ and $b$ are the two fixed points of $f$. By Theorem~\ref{T: two} (i) we can construct from $h$ an orientation reversing involution $\mu$ such that $\mu f^2\mu = f^{-2}$. We require that $\mu$, like $h$, interchanges $a$ and $b$; this is not given by the statement of Theorem~\ref{T: two}, however, it is immediate from the definition of $\mu$ in \eqref{E: invert} (where, in that equation, $s$ is our current homeomorphism $h$). Now define
\[
\tau(x)=
\begin{cases}
\mu(x) & \text{if $x\in[a,b]$},\\
f\mu f(x) & \text{if $x\in[b,a]$}.
\end{cases}
\] 
Then $\tau$ is an involution in \Hon and $\tau f\tau =f^{-1}$. 
\end{proof}

The hypothesis that $h$ interchanges the pair of fixed points of $f$  cannot be dropped from Theorem~\ref{T: two-}: one can construct examples of orientation reversing homeomorphisms $f$ and involutions $\tau$ such that $\tau f^2\tau=f^{-2}$ even though $f$ is not strongly reversible. There are also examples of orientation reversing circle homeomorphisms that are reversible, but not strongly reversible.


\bibliographystyle{plain}
\bibliography{reversible}

\end{document}